\documentclass[12pt]{amsart}
\usepackage{amsfonts,amssymb,amsmath,amsthm}
\usepackage{url}
\usepackage{enumerate}
\usepackage{bbm}
\usepackage{times}

\def\B{\mathcal{B}}
\def\D{\mathbb{D}}

\def\T{\partial \D}
\def\C{\mathbb{C}}
\def\Q{\mathcal{Q}}
\def\L{\mathcal L}
\def\qp{\Q_p}
\def\qk{\mathcal{Q}_K}
\def\f{\frac}
\begin{document}

\title
{ The pseudoanalytic extensions for some spaces of analytic functions}
\author{Guanlong Bao, Hasi Wulan and Fangqin Ye}
\address{
Department of Mathematics\\
    Shantou University\\
    Shantou, Guangdong,  515063,  China}
\email{glbaoah@163.com (G. Bao)}
\email{wulan@stu.edu.cn (H. Wulan)}
\email{yefqah@163.com (F. Ye)}

\thanks{The authors  are    supported  by NSF of China (No. 11371234).} 
\begin{abstract}
Using the Cauchy-Riemann operator, we characterize  $\qk$ spaces, Besov spaces and analytic Morrey spaces in terms of pseudoanalytic extensions of primitive functions. Our results are also true on  some classical Banach spaces, such as the Bloch space, $BMOA$ and the Dirichlet space.


\vskip 3mm
\end{abstract}

\maketitle

\section{Introduction}\vspace{0.2truecm}

Let $\D$ be the  unit disk in the complex plane $\C$.
Denote by  $H(\D)$  the space of functions analytic in $\D$. The Green function
in the unit disk  with singularity at $a\in\D$ is given by
$$
g(a, z) = \log\frac{1}{|\sigma_a(z)|}, \ z\in \D.
$$
Here
$$
\sigma_a(z)=\f{a-z}{1-\overline{a}z},
$$
is a M\"obius transformation of $\D$.

Throughout this paper, we
assume that $K:[0,\infty)\rightarrow [0, \infty)$ is a
right-continuous and increasing function. A function
$f\in H(\D)$ belongs to the space $\Q_K$ if
$$
\|f\|_{\qk}^2 = \sup_{a\in\D}\,\int_{\D}
|f'(z)|^{2}K\left(g(a, z)\right)dA(z)<\infty,
$$
where $dA(z)$  is the Lebesgue measure in $\D$.
By \cite[Theorem 2.1]{EW}, $f\in\qk$ if and only if
$$
\sup_{a\in\D}\int_\D |f'(z)|^2K\left(1-|\sigma_a(z)|^2\right)dA(z)<\infty.
$$
See \cite{EW} and
\cite{EWX} for more results of $\Q_K$ spaces. If
$K(t)=t^p$, $0\leq p<\infty$, then the space $\Q_K$ gives the space
$\Q_p$ (see \cite{AXZ, X1, X2}). In particular, $\Q_0$ is the Dirichlet space $\mathcal D$;
$\Q_1=BMOA$, the space of bounded mean oscillation (see \cite{Ba, Gi});  by \cite{AL}, for all $p\in (1, \infty)$,   the spaces $\qp$ are the same and equal to the Bloch space $\B$ which consists of all functions $f\in H(\D)$ with
$$
\|f\|_{\B}=\sup_{z\in \D}(1-|z|^2)|f'(z)|<\infty.
$$

For $0<p<\infty$, the Hardy space $H^p$ consists of all functions $f\in H(\D)$ with
$$
\|f\|_{H^p}^p=\sup_{0<r<1}\f{1}{2\pi}\int_0^{2\pi}|f(re^{i\theta})|^pd\theta<\infty.
$$
It is well-known that if $f\in H^p$, then its nontangential limit $f(e^{i\theta})$ exists almost everywhere (see \cite{D, G}).

In this paper, we need  the Cauchy-Riemann operator
$$
\overline{\partial}=\f{\partial}{\partial \overline{z}}=\f{1}{2}\left(\f{\partial}{\partial x}+i \f{\partial}{\partial y}\right), \ z=x+iy.
$$
 Dyn'kin \cite{Dyn} said that the so-called pseudoanalytic extension is a method of extension with a given estimate of the Cauchy-Riemann operator.
Dyn'kin \cite{Dyn} characterized many classical smoothness spaces applying  the pseudoanalytic extension method. Note that
$$
\bigcup_{0<p<1}\qp\subseteq BMOA\subseteq \bigcap_{0<q<\infty} H^q.
$$
Dyakonov and Girela \cite{Dy2} obtained  an interesting characterization of $\qp$ spaces in terms of pseudoanalytic extension  as follows.

\newtheorem*{thmA}{Theorem A}
\begin{thmA} If $0<p<1$ and $f\in \bigcap_{0<q<\infty} H^q$, then the following conditions are equivalent.
\begin{enumerate}
\item [(i)] $f\in\Q_p$.
\item[(ii)]$$\sup_{a\in\D}\int_\D |f'(z)|^2\left(\f{1}{|\sigma_a(z)|^2}-1\right)^pdA(z)<\infty.$$
\item[(iii)] There exists a function $\widetilde{f}\in C^1(\C\setminus\overline{\D})$ satisfying
$$
\widetilde{f}(z)=O(1), \ \text{as} \ \ z\rightarrow\infty,
$$
$$
\lim_{r\rightarrow1^+}\widetilde{f}(re^{i\theta})=f(e^{i\theta}), \ \text{a.e. and in} \ L^q([-\pi,\ \pi]) \ \text{for all} \ \ q\in[1,\infty),
 $$
and
$$
\sup_{a\in\D}\int_{\C\setminus\overline{\D}}|\overline{\partial}\widetilde{f}(z)|^2\left(|\sigma_a(z)|^2-1\right)^pdA(z)<\infty.
$$
\end{enumerate}
\end{thmA}
The proof of the above theorem in \cite{Dy2} involved the Calder\'on-Zygmund operators and the Muckenhoupt weights (see \cite{S}).
For  $a\in \D$ and $0<p<1$,  let
$$
U_a(z)=\left|1-\frac{1}{|\sigma_a(z)|^2}\right|^p=\frac{(1-|a|^2)^p||z|^2-1|^p}{|z-a|^{2p}}, \ z\in\C.
$$
 Dyakonov and Girela  \cite{Dy2} showed that $U_a$ is a
Muckenhoupt weight  with
$$
\sup_Q \left[\frac{1}{|Q|}\int_Q U_a(z)dA(z)\right]\left[\frac{1}{|Q|}\int_Q (U_a(z))^{-1}dA(z)\right]<\infty
$$
for all $a\in \D$. Here $Q$ ranges over the disks in $\C$ and $|Q|$ denotes the area of $Q$. The  above estimate of $U_a$ with $a=0$ was  also used to establish the $\qp$ corona theorem (see
\cite[Theorem 3.1]{X}). Applying Theorem A, Dyakonov and Girela \cite{Dy2} also gave some  nice properties of $\qp$ spaces.

Motivated by Theorem A, it is nature to ask some  questions as follows.

\vskip 1mm
\noindent{\bf Question 1.} Can we obtain pseudoanalytic extensions for all $\qp$, $0\leq p<\infty$?

\vskip 1mm
\noindent{\bf Question 2.} Let $n$ be a positive integer.  Does there exist a pseudoanalytic extension characterization of $\qp$ spaces such that the extension function is  $O(z^n)$ as $ z\rightarrow\infty$?

\vskip 1mm
\noindent{\bf Question 3.} Can we give characterizations for some spaces of analytic functions such as $\qk$ spaces, Besov spaces and analytic Morrey spaces in terms of pseudoanalytic extensions?

Throughout this article, denote by  $F$  the primitive function of $f\in H(\D)$; that is
$$
F(z)=\int_0^z f(w)dw, \  z\in \D.
$$
Related to the above questions, in this paper, we  characterize some analytic function spaces on $\D$ by  pseudoanalytic extension of primitive functions. Our method, without using
Calder\'on-Zygmund operators and Muckenhoupt weights, can be applied to some spaces, such as $\qk$ spaces, Besov spaces, analytic Morrey spaces and so on. Our results are also new for $\qp$ spaces.

In this paper, the symbol $A\approx B$ means that $A\lesssim B\lesssim
A$. We say that $A\lesssim B$ if there exists a constant $C$ such that $A\leq CB$.

\section{ $\qk$ spaces and pseudoanalytic extension }

In this section, we need two more constraints on $K$ as follows.
$$
\int_0^1\frac{\varphi_K(s)}{s}ds<\infty   \eqno{(2.1)}
$$
and
$$
\int_1^\infty\frac{\varphi_K(s)}{s^{1+q}}ds<\infty, \ \ \ \ 0<q<
2, \eqno{(2.2)}
$$
where
$$
\varphi_K(s)=\sup_{0\leq t\leq 1}\frac{K(st)}{K(t)},\ \ \ \ \ 0<s<\infty.
$$
Under conditions (2.1) and (2.2), $\qk$ spaces have been studied extensively (see \cite{WY, WZh1, WZh2}). From now on we always assume the function $K$ satisfying the double condition, namely $K(2t)\approx K(t)$ for all $t\in (0, 1)$.

A very useful tool in the study of $\qk$ spaces is $K$-Carleson
measure. Let $\ell(I)$ be the length of an arc $I$ of the unit circle $\T$.  Define the
Carleson box by
$$
S_G(I)=\begin{cases}
\{r\zeta\in G: 1-\f{\ell(I)}{2\pi}<r<1, \zeta\in I\}, & G=\D, \\
              \{r\zeta\in G: 1<r<1+\ell(I), \zeta\in I\}, & \enspace G=\C\setminus
\overline{\D}.  \end{cases}
$$
Following \cite{EWX} and \cite{WY}, a positive Borel measure $\mu$ on $G=\D$ or $G=\C\setminus
\overline{\D}$ is said to be a
$K$-Carleson measure if
$$
\sup_{I\subset \T}\int_{S_G(I)}K\left(\frac{|1-|z||}{\ell(I)}\right)d\mu(z)<\infty.
$$

By \cite{EW}, we know that all $\qk$ spaces are subsets of the Bloch space. Then the primitive function $F$ of a $\qk$ function $f$ must be in the Hardy space $H^2$ since the Bloch functions' Taylor coefficients is bounded (see \cite{An}). Therefore, the primitive function $F$ has its nontangential limit $F(e^{i\theta})$ almost everywhere on the unit circle.

The following is the main result of this section.

\newtheorem*{thm2.1}{Theorem 2.1}
\begin{thm2.1}   Suppose that $K$ satisfies (2.1) and (2.2). Let $f\in H(\D)$ with its primitive function $F\in H^2$ and let $n\geq 2$ be an integer.  Then the following conditions are equivalent.
\begin{enumerate}
\item [(i)] $f\in\qk$.
\item[(ii)] There exists a function $\widetilde{F_n}\in C^1(\C\setminus\overline{\D})$ satisfying
$$\lim_{r\rightarrow1^+}\widetilde{F_n}(re^{i\theta})=F(e^{i\theta})\  a.e.\  \theta \in [0, 2\pi], \eqno (2.3)$$
 $$\widetilde{F_n}(z)=O(z^n), \ \text{as} \ \ z\rightarrow\infty,  \eqno (2.4)$$
 $$\overline{\partial}\widetilde{F_n}(z)=O(z^{n-2}), \ \text{as} \ \ z\rightarrow\infty,  \eqno (2.5)$$
and
$$
\sup_{a\in\D}\int_{\C\setminus\overline{\D}}\f{|\overline{\partial}\widetilde{F_n}(z)|^2}{(|z|^n-1)^2}K\left(1-\frac{1}{|\sigma_a(z)|^2}\right)dA(z)<\infty.
\eqno (2.6)
$$
\item[(iii)] There exists a function $\widetilde{F_n}\in C^1(\C\setminus\overline{\D})$ satisfying (2.3),  (2.4), (2.5) and $|\overline{\partial}\widetilde{F_n}(z)|^2/(|z|^n-1)^2$ is a $K$-Carleson measure on  $\C\setminus
\overline{\D}$.
\end{enumerate}
\end{thm2.1}

Before embarking into the proof of Theorem 2.1, we state  some lemmas below.

\newtheorem*{lem2.2}{Lemma 2.2}
\begin{lem2.2} Suppose that  $K$ satisfies (2.1). Let $\mu$ be  a positive Borel measure on $\D$. Then the following conditions are equivalent.
\begin{enumerate}
\item [(i)] $\mu$ is a $K$-Carleson measure on  $\D$.
\item[(ii)]
$$
\sup_{a\in\D}\int_{\D} K(1-|\sigma_a(z)|^2)d\mu(z)<\infty.
$$
\item[(iii)]
$$
\sup_{a\in\D}\int_{\D} K\left(\frac{1}{|\sigma_a(z)|^2}-1\right)d\mu(z)<\infty.
$$
\end{enumerate}
\end{lem2.2}

\begin{proof}
Note that $(i)\Leftrightarrow(ii)$ was proved in \cite{EWX}. To show the equivalence of $(ii)$ and $(iii)$, let  $$
K_1(t)=K\left(\f{t}{1-t}\right), \ 0<t<\f{1}{2}.
$$
By the monotonicity of $K$,  $K_1(s)$ is also increasing for $s\in (0, 1)$.  Moreover,
$$
K(t)\leq K\left(\f{t}{1-t}\right) \leq K(2t), \ 0<t<\f{1}{2}.
$$
Since $K(2t)\approx K(t)$ for $t\in (0, 1)$, $K_1(t)\approx K(t)$ for $t\in (0, 1/2)$. Thus, the proof is complete. \end{proof}

Wulan and Zhu \cite{WZhu} proved that if $K$ satisfies (2.1) and  $n$ is a positive integer,  then $f\in \qk$ if and only if
$$
\sup_{a\in\D}\int_\D |f^{(n)}(z)|^2 (1-|z|^2)^{2n-2}K\left(1-|\sigma_a(z)|^2\right)dA(z)<\infty. \eqno(2.7)
$$
Combining this with  Lemma 2.2, we obtain immediately  a characterization of $\qk$ spaces as follows.  If $K$ satisfies (2.1) and $n$ is a positive integer, then  $f\in \qk$ if and only if
$$
\sup_{a\in\D}\int_\D |f^{(n)}(z)|^2 (1-|z|^2)^{2n-2} K\left(\f{1}{|\sigma_a(z)|^2}-1\right)dA(z)<\infty.
$$

For the case of $\mu$ on $\C\setminus \overline{\D}$, we give a similar description of $K$-Carleson measure  as follows.

\newtheorem*{lem2.3}{Lemma 2.3}
\begin{lem2.3} Suppose that  $K$ satisfies (2.1). Let $\mu$ be  a positive Borel measure on $\C\setminus \overline{\D}$. Then the following conditions are equivalent.
\begin{enumerate}
\item [(i)] $\mu$ is a $K$-Carleson measure on  $\C\setminus \overline{\D}$.
\item[(ii)]
$$
\sup_{a\in\D}\int_{1<|z|<1+2\pi} K\left(|\sigma_a(z)|^2-1\right)d\mu(z)<\infty.
$$
\item[(iii)]
$$
\sup_{a\in\D}\int_{1<|z|<1+2\pi} K\left(1-\frac{1}{|\sigma_a(z)|^2}\right)d\mu(z)<\infty.
$$
\end{enumerate}
\end{lem2.3}

\begin{proof}
$(i)\Rightarrow(ii)$.  Fix $a=re^{i\theta}\in\D$.  Let $I$ be the arc of center $e^{i\theta}$ and $\ell(I)=\frac{2\pi(1-|a|)}{2\pi |a|+1}$.  Then for any $z\in S_{\C\setminus \overline{\D}}(I)$, one gets
$$
|1-\bar{a}z|\geq 1-(1+\ell(I))|a|=\frac{1}{2\pi}\ell(I)
$$
and
\begin{eqnarray*}
|1-\bar{a}z|&\leq&|a|\left|e^{i\left(\theta+\frac{\ell(I)}{2}\right)}-\frac{1}{\bar{a}}\right|\\
&=&r\left|e^{i\frac{\ell(I)}{2}}-\frac{1}{r}\right|
=\left[(1-r)^2+4r\left(\sin\frac{\ell(I)}{4}\right)^2\right]^{\frac{1}{2}}\\
&\leq&\left[\left(\frac{2\pi r+1}{2\pi}\ell(I)\right)^2+4r\left(\frac{\ell(I)}{4}\right)^2\right]^{\frac{1}{2}}
\leq 2\ell(I).
\end{eqnarray*}
Hence
$$
|1-\bar{a}z|\approx 1-|a|\approx \ell(I),\ \ z\in S_{\C\setminus \overline{\D}}(I).
$$
Set
$$
S_n=\left\{r\zeta\in \C\setminus \overline{\D}: 1<r<1+\ell(2^nI), \zeta \in 2^nI\right\},
$$
where $2^nI$ denotes the arc with the same center as $I$ and with length $2^n\ell(I)$. Then
\begin{eqnarray*}
&~&\int_{1<|z|<1+2\pi} K\left(|\sigma_a(z)|^2-1\right)d\mu(z)\\
&\leq&\int_{S_{\C\setminus \overline{\D}}(I)}K\left(|\sigma_a(z)|^2-1\right)d\mu(z)
 +\sum_{n=1}^\infty \int_{S_n\backslash S_{n-1}}K\left(|\sigma_a(z)|^2-1\right)d\mu(z)\\
&\lesssim&\int_{S_{\C\setminus \overline{\D}}(I)}K\left(\frac{|z|-1}{\ell(I)}\right)d\mu(z)+\sum_{n=1}^\infty \int_{S_n\backslash S_{n-1}}K\left(\frac{|z|-1}{2^{2n}\ell(I)}\right)d\mu(z)\\
&\lesssim&\int_{S_{\C\setminus \overline{\D}}(I)}K\left(\frac{|z|-1}{\ell(I)}\right)d\mu(z)+\sum_{n=1}^\infty\sup_{z\in S_n} \frac{K\left(\frac{|z|-1}{2^{2n}\ell(I)}\right)} {K\left(\frac{|z|-1}{2^n\ell(I)}\right)}\int_{S_n}K\left(\frac{|z|-1}{\ell(2^nI)}\right)d\mu(z)\\
&\lesssim&\int_{S_{\C\setminus \overline{\D}}(I)}K\left(\frac{|z|-1}{\ell(I)}\right)d\mu(z)+\sum_{n=1}^\infty \varphi_K(2^{-n})\int_{S_n}K\left(\frac{|z|-1}{\ell(2^nI)}\right)d\mu(z).
\end{eqnarray*}
Since $d\mu$ is a $K$-Carleson measure and
$$
\sum_{n=1}^\infty \varphi_K(2^{-n})\thickapprox \int_0^1 \varphi_K(s)\frac{ds}{s}<\infty,
$$
 we obtain that
$$
\sup_{a\in\D}\int_{1<|z|<1+2\pi} K\left(|\sigma_a(z)|^2-1\right)d\mu(z)<\infty.
$$

$(ii)\Rightarrow(iii)$.  It follows from that
$$
K\left(1-\frac{1}{|\sigma_a(z)|^2}\right)\leq K\left(|\sigma_a(z)|^2-1\right)
$$
holds for all $a\in\D$ and $1<|z|<1+2\pi$.

$(iii)\Rightarrow(i)$.  For given a subarc $I\in \T$, set $e^{i\theta}$ the midpoint of $I$ and set
$$
a=\f{2\pi-\ell(I)}{2\pi(\ell(I)+1)}  e^{i\theta}.
$$
Then for any  $z\in S_{\C\setminus \overline{\D}}(I)$, we have
$$
0<|\sigma_a(z)|^2-1
\leq\frac{\frac{1}{\pi}(2\pi |a|+1)(2+\ell(I))\left(\ell(I)\right)^2}{\left(\frac{1}{2\pi}\ell(I)\right)^2}\leq8\pi(2\pi+1)(\pi+1).
$$
Let
$$
K_2(t)=K\left(\frac{t}{1+t}\right), \ 0<t<8\pi(2\pi+1)(\pi+1).
$$
By the monotonicity of $K$,  $K_2(s)$ is also increasing for $s\in (0, 1)$.  Moreover,
$$
K\left(\frac{t}{8\pi(2\pi+1)(\pi+1)+1}\right)\leq K_2(t) \leq K(t), \ 0<t<8\pi(2\pi+1)(\pi+1).
$$
Since  $K(2t)\approx K(t)$ for $t\in (0, 1)$,
$K_2(t)\approx K(t)$ for $t\in (0, 8\pi(2\pi+1)(\pi+1))$. Thus,
\begin{eqnarray*}
\int_{S_{\C\setminus \overline{\D}}(I)}K\left(\frac{|z|-1}{\ell(I)}\right)d\mu(z)
&\approx&\int_{S_{\C\setminus \overline{\D}}(I)}K\left(|\sigma_a(z)|^2-1\right)d\mu(z)\\
&\approx&\int_{S_{\C\setminus \overline{\D}}(I)}K\left(\frac{|\sigma_a(z)|^2-1}{|\sigma_a(z)|^2}\right)d\mu(z)\\
&\lesssim&\int_{1<|z|<1+2\pi} K\left(1-\frac{1}{|\sigma_a(z)|^2}\right)d\mu(z).
\end{eqnarray*}
Hence, (i) follows.
\end{proof}

From \cite{EWX} and \cite{WZh1}, if $K$ satisfies (2.1) and (2.2), then there exists a small enough constant $c$, $0<c<q/2$, such that $K(t)/t^c$ is increasing and $K(t)/t^{q-c}$ is decreasing in (0, 1). To prove Theorem 2.1, we also  need the following estimate.

\newtheorem*{lem2.4}{Lemma 2.4}
\begin{lem2.4} Let (2.1) and (2.2) hold for $K$. If  $s<\min (1+c, 2-q+c)$ and $2s+r-4\geq 0$, then
$$
\int_{\D} \frac{K\left(1-|\sigma_a(w)|^2\right)}{(1-|w|^2)^{s}|1-\overline{w}z|^r}dA(w)\lesssim
\frac{K\left(1-|\sigma_a(z)|^2\right)}{(1-|z|^2)^{s+r-2}}
$$
for all $a, z\in \D$. Here $c$ is a small enough positive  constant depending  only on  (2.1) and (2.2).
\end{lem2.4}
\begin{proof} For fixed $a$, $z\in \D$,  let $\lambda=\sigma_z(a)$. Then
$$
|\sigma_a(w)|=| \sigma_\lambda\circ \sigma_z(w)|.
$$
Note that  $2s+r-4\geq 0$. Checking the proof of Lemma 2.1 in \cite{BL}, one gets
\begin{eqnarray*}
&~&\int_{\D} \frac{K\left(1-|\sigma_a(w)|^2\right)}{(1-|w|^2)^{s}|1-\overline{w}z|^r}dA(w)\\
&\lesssim&\frac{K\left(1-|\lambda|^2\right)}{(1-|z|^2)^{s+r-2}}\int_{\D} \varphi_K\left(\frac{1-|u|^2} {|1-\overline{\lambda}u|^2}\right) \frac{1}{(1-|u|^2)^{s}}dA(u).
\end{eqnarray*}
Since $K$ satisfies (2.1) and (2.2), there exists a small enough positive  constant $c$ depending  only on  (2.1) and (2.2) (see  \cite{EWX, WZh1}),  such that
$$
\varphi_{K}(t)\lesssim t^c, \ 0<t\leq1
$$
and
$$
\varphi_{K}(t)\lesssim t^{q-c}, \ t\geq1.
$$
Thus,
\begin{eqnarray*}
&~&\int_{\D} \varphi_K\left(\frac{1-|u|^2} {|1-\overline{\lambda}u|^2}\right) \frac{1}{(1-|u|^2)^{s}}dA(u)\\
&\lesssim& \int_{\D}  \frac{(1-|u|^2)^{c-s}}{|1-\overline{\lambda}u|^{2c}}dA(u)
+\int_{\D}  \frac{(1-|u|^2)^{q-c-s}}{|1-\overline{\lambda}u|^{2q-2c}}dA(u).
\end{eqnarray*}
For $s<1+c<q+1-c$ and $s<2-q+c$, using Lemma 3.10  in \cite{Zh}, we get
$$
\int_{\D}  \frac{(1-|u|^2)^{c-s}}{|1-\overline{\lambda}u|^{2c}}dA(u)
+\int_{\D}  \frac{(1-|u|^2)^{q-c-s}}{|1-\overline{\lambda}u|^{2q-2c}}dA(u)\lesssim 1.
$$
Thus,
$$
\int_{\D} \frac{K\left(1-|\sigma_a(w)|^2\right)}{(1-|w|^2)^{s}|1-\overline{w}z|^r}dA(w)
\lesssim \frac{K\left(1-|\lambda|^2\right)}{(1-|z|^2)^{s+r-2}}\approx\frac{K\left(1-|\sigma_a(z)|^2\right)}{(1-|z|^2)^{s+r-2}}.
$$
\end{proof}

Now we are ready to finish the proof of Theorem 2.1.
\vskip 2mm

$(i)\Rightarrow (ii)$. From now on, we write
 $$
  w^*=1/\overline{w}, \ w\in \C\setminus\{0\}.
  $$
Let $f\in \qk$ and
$$
F(z)=\int_0^z f(\zeta)d\zeta, \  z\in \D.
$$
Set
$$
\widetilde{F_n}(z)=\sum_{i=0}^{n}\f{(-1)^i}{i!}(z^*-z)^iF^{(i)}(z^*), \ z\in \C\setminus\overline{\D}.
$$
Clearly, $\widetilde{F_n}\in C^1(\C\setminus\overline{\D})$ satisfying
$$
\lim_{r\rightarrow1^+}\widetilde{F_n}(re^{i\theta})=F(e^{i\theta}) \ \ a.e.\ \theta\in [0, 2\pi],
$$
and
 $$
 \widetilde{F_n}(z)=O(z^n), \ \text{as} \ \ z\rightarrow\infty.
 $$
Note that
$$
\overline{\partial}\widetilde{F_n}(z)=\f{(-1)^{n+1}}{n!}(z^*-z)^n(z^*)^2F^{(n+1)}(z^*).
$$
Then
$$
\overline{\partial}\widetilde{F_n}(z)=O(z^{n-2}), \ \text{as} \ \ z\rightarrow\infty.
$$
Making  the change of variable $z=\zeta^*$, we deduce that
\begin{eqnarray*}
&~&\sup_{a\in\D}\int_{\C\setminus\overline{\D}}\f{\left|\overline{\partial}\widetilde{F_n}(z)\right|^2}{(|z|^n-1)^2}K\left(1-\frac{1}{|\sigma_a(z)|^2}\right)dA(z)\\
&=&\f{1}{(n!)^2}\sup_{a\in\D}\int_{\C\setminus\overline{\D}}\f{|z^*-z|^{2n}|z^*|^4\left|f^{(n)}(z^*)\right|^2}{(|z|^n-1)^2}K\left(1-\frac{1}{|\sigma_a(z)|^2}\right)dA(z)\\
&\approx&\sup_{a\in\D}\int_{\D} \left|f^{(n)}(\zeta)\right|^2(1-|\zeta|)^{2n-2}K\left(1-|\sigma_a(\zeta)|^2\right)dA(\zeta).
\end{eqnarray*}
Combining this with  (2.7), we get the  desired result.

$(ii)\Rightarrow (i)$.  Assume that (ii) holds. Let $z\in \D$ and $R>1$.  Following a technique in \cite{Dy1} or \cite{Dy2} and in view of (2.3), we employ the Cauchy-Green formula to the function that equals $F$ in $\D$ and $\widetilde{F_n}$ in $\C\setminus\overline{\D}$. Then
$$
 F(z)=\f{1}{2\pi i} \int_{|\zeta|=1}\f{F(\zeta)}{\zeta-z}d\zeta=\f{1}{2\pi i} \int_{|\zeta|=R}\f{\widetilde{F_n}(\zeta)}{\zeta-z}d\zeta-\f{1}{\pi}\int_{1<|\zeta|<R} \f{\overline{\partial}\widetilde{F_n}(\zeta)}{\zeta-z}dA(\zeta).   \eqno{(2.8)}
$$
By (2.4),
 $$
 \int_{|\zeta|=R}\f{\widetilde{F_n}(\zeta)}{(\zeta-z)^{n+2}}d\zeta\rightarrow0, \ \text{as} \ R\rightarrow \infty.
 $$
 By (2.5), we know
 $$
 \left|\int_{\C\setminus\overline{\D}}\f{\overline{\partial}\widetilde{F_n}(\zeta)}{(\zeta-z)^{n+2}}dA(\zeta)\right|<\infty.
 $$
 These together with  (2.8) give
 $$
 F^{(n+1)}(z)=-\f{(n+1)!}{\pi}\int_{\C\setminus\overline{\D}}\f{\overline{\partial}\widetilde{F_n}(\zeta)}{(\zeta-z)^{n+2}}dA(\zeta).
 $$
By the H\"older inequality, the change of variable and the well known estimates in Zhu's book \cite[Lemma 3.10]{Zh}, we obtain
\begin{eqnarray*}
\left|F^{(n+1)}(z)\right|^2&\lesssim & \int_{\C\setminus\overline{\D}}\f{1}{|\zeta-z|^4}dA(\zeta)\int_{\C\setminus\overline{\D}}\f{|\overline{\partial}\widetilde{F_n}(\zeta)|^2}{|\zeta-z|^{2n}}dA(\zeta)\\
&\approx& \int_{\D} \f{1}{|1-\overline{w}z|^4}dA(w)\int_{\C\setminus\overline{\D}}\f{|\overline{\partial}\widetilde{F_n}(\zeta)|^2}{|\zeta-z|^{2n}}dA(\zeta)\\
&\approx& \f{1}{(1-|z|^2)^2}\int_{\C\setminus\overline{\D}}\f{|\overline{\partial}\widetilde{F_n}(\zeta)|^2}{|\zeta-z|^{2n}}dA(\zeta).
\end{eqnarray*}
Using  Lemma 2.4,  $n\geq 2$ and (2.6), we see that
\begin{eqnarray*}
&~&\sup_{a\in\D}\int_{\D}(1-|z|^2)^{2n-2}\left|F^{(n+1)}(z)\right|^2K\left(1-|\sigma_a(z)|^2\right)dA(z)\\
&\lesssim&\sup_{a\in\D}\int_{\D}\int_{\C\setminus\overline{\D}}\f{|\overline{\partial}\widetilde{F_n}(\zeta)|^2}{|\zeta-z|^{2n}}dA(\zeta)(1-|z|^2)^{2n-4}K\left(1-|\sigma_a(z)|^2\right)dA(z)\\
&\lesssim&\sup_{a\in\D}\int_{\D}\int_{\D}\f{K\left(1-|\sigma_a(z)|^2\right)}{|1-\overline{w}z|^4}dA(z)|\overline{\partial}\widetilde{F_n}(w^*)|^2|w|^{2n-4}dA(w)\\
&\lesssim&\sup_{a\in\D}\int_{\D}\frac{|\overline{\partial}\widetilde{F_n}(w^*)|^2|w|^{2n-4}}{(1-|w|^2)^2}K\left(1-|\sigma_a(w)|^2\right)dA(w)\\
&\approx&\sup_{a\in\D}\int_{\C\setminus\overline{\D}}\frac{|\overline{\partial}\widetilde{F_n}(\zeta)|^2}{(|\zeta|^n-1)^2}K\left(1-\frac{1}{|\sigma_a(\zeta)|^2}\right)dA(\zeta)<\infty.
\end{eqnarray*}
Note that $F'(z)=f(z)$ for $z\in \D$. By (2.7),  we get $f\in \qk$.

$(ii)\Leftrightarrow (iii)$. If
$$
\overline{\partial}\widetilde{F_n}(z)=O(z^{n-2}), \ \text{as} \ \ z\rightarrow\infty,
$$
then
$$
\sup_{a\in\D}\int_{|z|\geq 1+2\pi}\f{|\overline{\partial}\widetilde{F_n}(z)|^2}{(|z|^n-1)^2}K\left(1-\frac{1}{|\sigma_a(z)|^2}\right)dA(z)<\infty.
$$
This together with Lemma 2.3, we see that condition  (ii) is equivalent to condition (iii).
We finish the proof.
\hfill{$\square$}

\vskip 2mm
\noindent{\bf Remark.} In Theorem 2.1, by the particular choice of $K$, we obtain the corresponding characterizations of $\qp$ for $0<p<1$, $BMOA$ and the Bloch space.

\section {Besov spaces, Morrey spaces and pseudoanalytic extension}

In this section, we describe Besov spaces and Morrey spaces by pseudoanalytic extension of primitive functions.

For $1<p<\infty$, the Besov space $B_p$ is the space of analytic functions $f$ in $\D$ such that
$$
\|f\|_{B_p}^p=\int_{\D}|f'(z)|^p(1-|z|^2)^{p-2}dA(z)<\infty.  \eqno(3.1)
$$
In particular, if $p=2$, then $B_p$ is the Dirichlet space $\mathcal D$.
Note that  (3.1) does not hold for $p=1$.   The Besov space $B_1$ consists of all
analytic functions $f$ on $\D$ which have a representation as
$$
f(z)=\sum_{k=1}^{\infty}c_k\sigma_{a_k}(z), \,\, a_k\in \D \,\,\,\,
\text{and} \,\,\sum_{k=1}^{\infty}|c_k|<\infty.
$$
Form \cite{AFP}, $f\in B_1$ if and only if
$$
\int_\D |f''(z)|dA(z)<\infty.
$$
Let  $1\leq p<\infty$ and $n\geq2$ be an integer. It is well known that $f\in B_p$ if and only if
$$
\int_{\D}|f^{(n)}(z)|^p(1-|z|^2)^{np-2}dA(z)<\infty. \eqno(3.2)
$$
See \cite{Da, Zh1, Zh} for more results of Besov spaces.

Note that all Besov spaces are subsets of the Bloch space. Thus, if $f\in B_p$, then its primitive function $F$ belongs to the Hardy space $H^2$.
Now we give a pseudoanalytic extension  characterization of Besov spaces as follows.

\newtheorem*{thm3.1}{Theorem 3.1}
\begin{thm3.1}   Let $f\in H(\D)$ with its primitive function $F\in H^2$ and let $n\geq 2$ be an integer. If $1\leq p <\infty$, then the following conditions are equivalent.
\begin{enumerate}
\item [(i)] $f\in B_p$.
\item[(ii)] There exists a function $\widetilde{F_n}\in C^1(\C\setminus\overline{\D})$ satisfying
$$\lim_{r\rightarrow1^+}\widetilde{F_n}(re^{i\theta})=F(e^{i\theta})\  a.e.\  \theta \in [0, 2\pi], \eqno (3.3)$$
 $$\widetilde{F_n}(z)=O(z^n), \ \text{as} \ \ z\rightarrow\infty,  \eqno (3.4)$$
 $$\overline{\partial}\widetilde{F_n}(z)=O(z^{n-2}), \ \text{as} \ \ z\rightarrow\infty,  \eqno (3.5)$$
and
$$
\int_{\C\setminus\overline{\D}}\f{|\overline{\partial}\widetilde{F_n}(z)|^p}{\left(|z|^{\frac{np}{2}-p+2}-1\right)^2}dA(z)<\infty.
\eqno (3.6)
$$
\end{enumerate}
\end{thm3.1}
\begin{proof}$(i)\Rightarrow (ii)$. Let $f\in B_p$.  Set
$$
\widetilde{F_n}(z)=\sum_{i=0}^{n}\f{(-1)^i}{i!}(z^*-z)^iF^{(i)}(z^*), \ z\in \C\setminus\overline{\D}.
$$
It is easy to check that $\widetilde{F_n}\in C^1(\C\setminus\overline{\D})$ and $\widetilde{F_n}$ satisfies  (3.3), (3.4) and (3.5). Furthermore, by the following estimates
\begin{eqnarray*}
&&\int_{\C\setminus\overline{\D}}\f{|\overline{\partial}\widetilde{F_n}(z)|^p}{\left(|z|^{\frac{np}{2}-p+2}-1\right)^2}dA(z)\\
&=&\f{1}{(n!)^p}\int_{\C\setminus\overline{\D}}\f{|z^*-z|^{np}|z^*|^{2p}\left|f^{(n)}(z^*)\right|^p}{\left(|z|^{\frac{np}{2}-p+2}-1\right)^2}dA(z)\\
&\approx&\int_{\D}\f{(1-|w|^2)^{np}\left|f^{(n)}(w)\right|^p}{\left(1-|w|^{\frac{np}{2}-p+2}\right)^2}dA(w)\\
&\approx&\int_{\D} (1-|w|)^{np-2}|f^{(n)}(w)|^pdA(w)
\end{eqnarray*}
and (3.2), we know that $\widetilde{F_n}$ also satisfies (3.6).

 $(ii)\Rightarrow (i)$.  Checking the proof of Theorem 2.1, one gets
 $$
 F^{(n+1)}(z)=-\f{(n+1)!}{\pi}\int_{\C\setminus\overline{\D}}\f{\overline{\partial}\widetilde{F_n}(\zeta)}{(\zeta-z)^{n+2}}dA(\zeta),\, z\in \D.
 $$
Applying  the H\"older inequality and  \cite[Lemma 3.10]{Zh}, we obtain
\begin{eqnarray*}
\left|F^{(n+1)}(z)\right|^p&\lesssim & \left(\int_{\C\setminus\overline{\D}}\f{1}{|\zeta-z|^4}dA(\zeta)\right)^{p-1} \int_{\C\setminus\overline{\D}}\f{|\overline{\partial}\widetilde{F_n}(\zeta)|^p}{|\zeta-z|^{np-2p+4}}dA(\zeta)\\
&\approx& \left(\int_{\D} \f{1}{|1-\overline{w}z|^4}dA(w)\right)^{p-1} \int_{\C\setminus\overline{\D}}\f{|\overline{\partial}\widetilde{F_n}(\zeta)|^p}{|\zeta-z|^{np-2p+4}}dA(\zeta)\\
&\approx& \f{1}{(1-|z|^2)^{2p-2}}\int_{\C\setminus\overline{\D}}\f{|\overline{\partial}\widetilde{F_n}(\zeta)|^p}{|\zeta-z|^{np-2p+4}}dA(\zeta).
\end{eqnarray*}
Note that  $n\geq 2$. Using \cite[Lemma 3.10]{Zh} again, we deduce that
\begin{eqnarray*}
&~&\int_{\D}(1-|z|^2)^{np-2}\left|F^{(n+1)}(z)\right|^pdA(z)\\
&\lesssim&\int_{\D}\int_{\C\setminus\overline{\D}}\f{|\overline{\partial}\widetilde{F_n}(\zeta)|^p}{|\zeta-z|^{np-2p+4}}dA(\zeta)(1-|z|^2)^{np-2p} dA(z)\\
&\approx&\int_{\D}\int_{\D}\f{|\overline{\partial}\widetilde{F_n}(w^*)|^p|w|^{np-2p}}{|1-\overline{w}z|^{np-2p+4}}dA(w)(1-|z|^2)^{np-2p} dA(z)\\
&\approx&\int_{\D}\int_{\D}\f{(1-|z|^2)^{np-2p}}{|1-\overline{w}z|^{np-2p+4}}dA(z)|\overline{\partial}\widetilde{F_n}(w^*)|^p|w|^{np-2p}dA(w)\\
&\lesssim&\int_{\D}\left(1-|w|^2\right)^{-2}|\overline{\partial}\widetilde{F_n}(w^*)|^p|w|^{np-2p}dA(w)\\
&\approx&\int_{\D}\left(1-|w|^{\frac{np}{2}-p+2}\right)^{-2}|\overline{\partial}\widetilde{F_n}(w^*)|^p|w|^{np-2p}dA(w)\\
&\approx&\int_{\C\setminus\overline{\D}}\frac{|\overline{\partial}\widetilde{F_n}(\zeta)|^p}{(|\zeta|^{\frac{np}{2}-p+2}-1)^2} dA(\zeta).
\end{eqnarray*}
Since  $F'(z)=f(z)$ for $z\in \D$,    by (3.2) and (3.6),  we obtain  $f\in B_p$. This finishes the proof.
\end{proof}

 For $0< \lambda \leq1$, the analytic Morrey space $\L^{2, \lambda}$ is the set of all functions $f\in H^2$ satisfying
$$
\|f\|_{\L^{2, \lambda}}=\sup_{I\subseteq \T} \left(\f{1}{\ell(I)^\lambda}\int_I|f(\zeta)-f_I|^2|d\zeta|\right)^{1/2}<\infty,
$$
where $\ell(I)$ denotes the length of $I$ and
$$
f_I=\f{1}{\ell(I)}\int_I f(\zeta)|d\zeta|.
$$
Clearly, $\L^{2, 1}$ is the space $BMOA$. By \cite[p. 54]{X2}, $f\in \L^{2, \lambda}$ if and only if
$$
\sup_{a\in \D} (1-|a|^2)^{1-\lambda}\int_\D |f'(z)|^{2}(1-|\sigma_a(z)|^2) dA(z)<\infty.
$$
 Let $n\geq1$ be an integer.        Using a similar statement  in \cite{AN}, we see that $f\in \L^{2, \lambda}$ if and only if
$$
\sup_{a\in \D} (1-|a|^2)^{1-\lambda}\int_\D |f^{(n)}(z)|^{2} (1-|z|^2)^{2n-2}(1-|\sigma_a(z)|^2) dA(z)<\infty. \eqno(3.7)
$$
Recently, the interest in analytic Morrey  space has grown rapidly. See \cite{C, LL, WX1, WX2, WX, XX, XY} for more results of $\L^{2, \lambda}$ spaces.

Note that if $f\in \L^{2, \lambda}$, then its primitive function $F$ belongs to $H^2$. Applying  (3.7), we get the following characterization of analytic Morrey space. The proof is similar to Theorem 2.1, we omit it.

\newtheorem*{thm3.2}{Theorem 3.2}
\begin{thm3.2}   Let $f\in H(\D)$ with its primitive function $F\in H^2$ and let $n\geq 2$ be an integer. If $0< \lambda \leq1$, then the following conditions are equivalent.
\begin{enumerate}
\item [(i)] $f\in\L^{2, \lambda}$.
\item[(ii)] There exists a function $\widetilde{F_n}\in C^1(\C\setminus\overline{\D})$ satisfying
$$\lim_{r\rightarrow1^+}\widetilde{F_n}(re^{i\theta})=F(e^{i\theta})\  a.e.\  \theta \in [0, 2\pi], $$
 $$\widetilde{F_n}(z)=O(z^n), \ \text{as} \ \ z\rightarrow\infty,  $$
 $$\overline{\partial}\widetilde{F_n}(z)=O(z^{n-2}), \ \text{as} \ \ z\rightarrow\infty,$$
and
$$
\sup_{a\in\D}(1-|a|^2)^{1-\lambda}\int_{\C\setminus\overline{\D}}\f{|\overline{\partial}\widetilde{F_n}(z)|^2}{(|z|^n-1)^2}
\left(1-\frac{1}{|\sigma_a(z)|^2}\right)dA(z)<\infty.
$$
\end{enumerate}
\end{thm3.2}

\vskip 2mm
{\bf  Remark.} When the paper is complete, we found that Wei and Shen  obtained the special case of our Theorem 2.1 for $K(t)=t$ and $n=2$. However, our proof in this paper and those of  \cite{WS} are different.

\end{document}